\documentclass[11pt]{amsart}

\usepackage{mathrsfs}
\usepackage{amscd}
\usepackage{amsmath}
\usepackage{amssymb}
\usepackage{amsthm}
\usepackage{epsf}
\usepackage{verbatim}
\usepackage{enumitem}
\usepackage{mathtools}
\usepackage{xcolor}
\usepackage{array}
\usepackage{bm}
\usepackage{bbm}
\usepackage{graphicx}
\input epsf.tex
\usepackage{amssymb}

\usepackage[left=3cm, right=3cm, top=2cm]{geometry}

\newtheorem{theorem}{Theorem}[section]
\newtheorem{lemma}[theorem]{Lemma}
\newtheorem{corollary}[theorem]{Corollary}
\theoremstyle{definition}

\theoremstyle{definition}

\theoremstyle{definition}

\theoremstyle{definition}
\newtheorem{definition}[theorem]{Definition}

\newcommand{\Nat}{\mathbb{N}}

\newcommand{\C}{\mathcal{C}}

\newcommand{\ve}{\varepsilon}

\newcommand{\hull}{\textup{hull}}
\newcommand{\fin}{\textup{fin}}
\newcommand{\odd}{\textup{odd}}
\newcommand{\Tutte}{\textup{Tutte}}
\newcommand{\Hall}{\textup{Hall}}

\title{Baire measurable matchings in non-amenable graphs}

\author{Alexander Kastner}
\address{UCLA Department of Mathematics}
\email{akastner@math.ucla.edu}

\author{Clark Lyons}
\address{E\"{o}tv\"{o}s Lor\'{a}nd University, Institute of Mathematics}
\email{lyons.clark@renyi.hu}

\begin{document}

\maketitle

\begin{abstract}
We prove that every Schreier graph of a free Borel action of a finitely generated non-amenable group admits a Baire measurable perfect matching, and that the Schreier graph of a free computable action of a finitely generated non-amenable group admits a computable perfect matching. These results were previously only known in the bipartite setting. To prove them, we establish variants of Tutte's theorem on perfect matchings in non-bipartite graphs. We also prove that every Borel non-amenable bounded-degree graph with only even degrees admits a Baire measurable balanced orientation.
\end{abstract}

\section{Introduction}

In this paper we consider non-amenable Borel graphs on Polish spaces and prove that certain classical combinatorial problems can be solved Baire measurably, that is, on a Borel comeager invariant set. Our method also gives analogous results for non-amenable highly computable graphs. Even in graphs on Polish spaces, our notion of non-amenability is entirely combinatorial, as explained in the following section. Our main theorem concerns the existence of perfect matchings:
\begin{theorem}\label{main}
Let $G$ be a Borel graph such that each component is an infinite, bounded-degree, non-amenable, vertex-transitive graph. Then $G$ admits a Borel perfect matching on a Borel comeager invariant set.
\end{theorem}

\begin{corollary}\label{maincor}
Every Schreier graph of a free Borel action of a finitely generated non-amenable group admits a Borel perfect matching on a Borel comeager invariant set.
\end{corollary}
In \cite{MUMatchings}, Marks and Unger studied Baire measurable matchings in the context of \emph{bipartite} Borel graphs, with a view toward applications to Baire measurable equidecompositions. Though not explicitly stated in their paper, their Theorem~1.3 implies that every bipartite Borel graph whose components are of bounded degree, regular, and non-amenable has a Baire measurable perfect matching. Thus, our theorem can be viewed as an extension of their result to the non-bipartite setting. The existence of regular, quasi-transitive, non-amenable graphs without perfect matchings (see Remark~22 of \cite{BHTOrientations}) leads us to assume vertex-transitivity, not just regularity.

Our method also yields analogous results for computable graphs.

\begin{theorem}\label{maincomp}
Let $G$ be a computable graph that is of bounded degree, vertex-transitive, and non-amenable. Then $G$ admits a computable perfect matching.
\end{theorem}

\begin{corollary}\label{compcor}
Every Schreier graph of a free computable action of a finitely generated non-amenable group admits a computable perfect matching.
\end{corollary}

These results all rely on technical variants of Tutte's theorem, which we state in Section~3. Techniques that produce graph labelings in both the Baire measurable and the highly computable case were studied by Qian and Weilacher in \cite{complab}. However, the fact that our method applies to both the Baire measurable and highly computable settings does not appear to follow from any previously known precise connections between these two settings.

Another theorem we prove concerns balanced orientations. Given a graph with only even degrees, a \textbf{balanced orientation} is an orientation of the edges such that every vertex has in-degree equal to out-degree. Euler's classical theorem on Euler circuits, together with a compactness argument, shows that every locally finite graph with only even degrees admits a balanced orientation. For non-amenable, bounded-degree Borel graphs, we have the following:
\begin{theorem}\label{orientation}
Let $G$ be a bounded-degree, non-amenable Borel graph with only even degrees. Then $G$ admits a Borel balanced orientation on a Borel comeager invariant set.
\end{theorem}

Our arguments draw inspiration from the study of factor-of-i.i.d.\ combinatorial structures for Cayley graphs of non-amenable groups, or more generally for quasi-transitive unimodular non-amenable graphs; see \cite{LN}, \cite{CL}, \cite{BHTOrientations}. By studying the spectrum of the Markov operator associated with random walks on these graphs, one shows a measure-expansion property for the associated Bernoulli graphings. This measure-expansion property is then used to establish the existence of a combinatorial structure (say, a perfect matching or balanced orientation) for the Bernoulli graphing on a Borel conull invariant set. While our arguments generally differ from those used to obtain factors of i.i.d., several of the ideas we use were inspired by that setting. It seems likely that many of the results pertaining to factors of i.i.d. for non-amenable graphs will have Baire measurable analogues.

\section{Notation}

We will use many standard notions from graph theory, recalled here, as well as some non-standard notions needed to formalize our technical variants of Tutte's theorem introduced in the next section. A \textbf{graph} $G$ will always be assumed to be symmetric and irreflexive. The \textbf{vertex set} of $G$ will be denoted by $V(G)$, and its \textbf{edge set} will be denoted by $E(G)$.

If $X$ and $Y$ are (disjoint) sets of vertices in a graph $G$, we let $$E(X,Y)=E(G)\cap(X\times Y)$$ denote the set of edges between $X$ and $Y$. We say that $G$ is \textbf{bipartite} if $V(G)$ can be partitioned into two sets $V(G)=X\sqcup Y$ such that $E(G)\subseteq E(X,Y)$. A \textbf{matching} is a set $M\subseteq E(G)$ such that each vertex is incident to at most one edge of $M$. We say $M$ is a \textbf{perfect matching} if each vertex is incident to exactly one edge of $M$. This paper concerns perfect matchings, and to emphasize that a matching is not assumed to be perfect we may use the term ``partial matching''.

If $X$ is a set of vertices, we let $G-X$ denote the graph with vertex set $$V(G-X)=V(G)-X$$ and edge set $$E(G-X)=E(G)\cap\Big(V(G-X)\times V(G-X)\Big).$$ We will also abuse notation and write $G-M$ for the graph $G-X$, where $X$ is the set of vertices incident to an edge of $M$.

All graphs considered in this paper will be \textbf{locally finite}, meaning that each vertex is related to only finitely many other vertices by the edge relation; that is, it has finitely many neighbors. A graph $G$ is \textbf{vertex-transitive} if any vertex can be mapped to any other vertex by a graph isomorphism from $G$ to $G$.

Given a group $G$ with a finite symmetric set of generators $S$ acting on a set $X$, we will consider the \textbf{Schreier graph} $G$ with vertex set $V(G)=X$ and edge set $$E(G)=\{(x,y)\in V(G)\times V(G):x\neq y\text{ and }\exists s\in S\text{ such that }s\cdot x=y\}.$$ Through this construction, free actions of finitely generated groups give rise to locally finite vertex-transitive graphs.

We say an infinite connected graph $G$ of bounded degree is \textbf{non-amenable} if there exists $\delta > 0$ such that whenever $F \subseteq V(G)$ is finite, the set of edges $E(F, V(G) \setminus F)$ between $F$ and $V(G) \setminus F$ satisfies $|E(F, V(G) \setminus F)| \ge \delta\,|F|$. For example, the Cayley graphs of finitely generated non-amenable groups with respect to any finite symmetric generating set (not containing the identity) are non-amenable graphs. There are many other examples that are not Cayley graphs. For example, non-amenable quasi-transitive unimodular graphs have been studied in \cite{LP} and \cite{BHTOrientations}. Some graphs that are not unimodular, such as the grandparent graph, are also covered by this definition.

\section{Perfect matchings}

The classical theorem of Tutte, repeated below, characterizes when a locally finite graph admits a perfect matching.
\begin{theorem}[Tutte] \label{TutteThm}
A locally finite graph $G$ admits a perfect matching if and only if whenever $X \subseteq V(G)$ is finite, the graph $G - X$ has at most $|X|$ finite components of odd size.
\end{theorem}
By \textbf{Tutte's condition} we will mean the condition that ``$G-X$ has at most $|X|$ odd components for each finite $X \subseteq V(G)$."

The proof of Theorem~\ref{main} consists of two steps. First, we establish a Baire measurable variant of Tutte's theorem which gives a sufficient condition for a locally finite Borel graph to admit a perfect matching on a Borel comeager invariant set (Theorem~\ref{BMTutte}). Second, we show that non-amenable, vertex-transitive graphs satisfy this sufficient condition (Lemma~\ref{NonamenableExpansion}).

\begin{definition}
If $G$ is a locally finite graph and $X \subseteq V(G)$ is finite, define
\[ \C_\fin(X) := \{\textup{finite components of $G - X$}\} \]
and
\[ \C_\odd(X) := \{\textup{odd components of $G - X$}\}.\]
Also let
\[ \hull_\fin(X) := X \cup \bigcup \C_\fin(X),\]
and
\[ \hull_\odd(X) := X \cup \bigcup \C_\odd(X).\]
We sometimes add superscripts to indicate the ambient graph when there is ambiguity.
\end{definition}


\begin{theorem} \label{BMTutte}
Let $G$ be a locally finite Borel graph on a Polish space $V(G)$, and suppose there exists $\ve > 0$ such that for every finite set $X \subseteq V(G)$, we have
\[ |X| \geq |\C_\odd(X)| + \ve \,|\hull_\odd(X)|.\]
Then $G$ admits a Borel perfect matching on a Borel comeager invariant set.
\end{theorem}

Note that when verifying Tutte's condition or this $\varepsilon$-strengthened variant in a locally finite graph $G$, it is sufficient to consider only finite sets of vertices $X$ such that $\hull_\odd(X)$ is connected (assuming that $G$ has no odd components). This is because the odd hull of a set can be decomposed into connected components, and Tutte's condition is additive on these components.

For the proof, we introduce a further refinement of Tutte's condition. A locally finite graph $G$ satisfies $\Tutte_{\varepsilon, k}$ if (i) Tutte's condition holds, and (ii) whenever $X \subseteq V(G)$ is finite such that $\hull_\odd(X)$ is connected with size at least $k$,
\[ |X| \geq |\C_\odd(X)| + \varepsilon \,|\hull_\odd(X)|.\]
Observe that the condition in Theorem~\ref{BMTutte} is equivalent to $\Tutte_{\ve, 1}$. This is an analogue of $\Hall_{\ve, k}$ in the proof of Theorem~1.3 in \cite{MUMatchings}. Our proof of Theorem~\ref{BMTutte} follows the same general strategy as the proof in \cite{MUMatchings}. In particular, we will need the following lemma from that paper.

\begin{lemma}[Marks-Unger]\label{sparselemma}
Let $G$ be a locally finite Borel graph on a Polish space $V(G)$, and let $f: \Nat \to \Nat$. Then there exist Borel sets $A_n \subseteq V(G)$, $n \in \Nat$, such that $\bigcup_n A_n$ is a Borel comeager invariant set and $d_G(x, y) > f(n)$ whenever $x, y$ are distinct vertices in $A_n$.
\end{lemma}

This decomposition into sparse Borel sets, combined with the following lemma, allow for the construction of a matching on a comeager invariant set.

\begin{lemma} \label{indstep}
Let $G$ be a locally finite graph satisfying $\textup{Tutte}_{\varepsilon,k}$. Let $\varepsilon'$ and $k'$ be such that
\begin{itemize}
    \item $0<\varepsilon'<\varepsilon$
    \item $k'\geq\max\{k,\frac{2}{\varepsilon-\varepsilon'}\}$
\end{itemize}
Let $A\subseteq V(G)$ be such that $d_G(x,y)>2k'$ for distinct $x,y\in A$. And for each $x\in A$ let $m(x)\in V(G)$ be a neighbor of $x$ such that the edge $e_x$ from $x$ to $m(x)$ belongs to a matching of $G$ covering the ball $B_{k'}(x)$.

If $M$ denotes the (partial) matching $$M=\{e_x:x\in A\},$$ then the graph $G-M$ satisfies $\Tutte_{\varepsilon',k'}$.

\end{lemma}

\begin{proof}
Assume first that $G-M$ has no odd components. This will be verified at the end of the proof.

Let $X\subseteq V(G-M)$ be finite with connected odd hull. Note that each $C\in\mathcal{C}_\odd^{G-M}(X)$ is not a component of $G-M$ and so must be adjacent to a vertex in $X$. It follows that $\mathcal{C}_\odd^{G-M}(X)$ and $\hull_\odd^{G-M}(X)$ are finite. List $\{x_1,\ldots,x_r\}$ as the set of vertices $x\in A$ such that either $x$ or $m(x)$ is adjacent (in $G$) to a component $C\in\mathcal{C}_\odd^{G-M}(X)$. Consider $$Y=X\cup\{x_1,m(x_1),\ldots,x_r,m(x_r)\}\subseteq V(G).$$ We have that $\mathcal{C}_\odd^G(Y)\supseteq\mathcal{C}_\odd^{G-M}(X)$ and $\hull_\odd^G(Y)$ is connected.

If $r=0$ then, for $X$, Tutte's condition and its $\varepsilon'$-strengthening in $G-M$ follow from Tutte's condition and its $\varepsilon$-strengthening in $G$. This is because $\varepsilon'<\varepsilon$ and $k'\geq k$.

If $r\geq2$ then, because elements of $A$ are at distance $>2k'$ and $\hull_\odd^G(Y)$ is connected, we have $|\hull_\odd^G(Y)|\geq rk'$. It follows that $|\hull_\odd^G(Y)|\geq k'$ and we can write $$|Y|\geq|\mathcal{C}_\odd^G(Y)|+\varepsilon\,|\hull_\odd^G(Y)|.$$ Now we can compute
\[
\begin{aligned}
|X| &= |Y| - 2r \\ 
&\geq|\mathcal{C}_\odd^G(Y)|+\varepsilon\,|\hull_\odd^G(Y)|-2r \\
&\geq|\mathcal{C}_\odd^{G-M}(X)|+\varepsilon'\,|\hull_\odd^{G-M}(X)|+(\varepsilon-\varepsilon')\,|\hull_\odd^G(Y)|-2r \\
&\geq|\mathcal{C}_\odd^{G-M}(X)|+\varepsilon'\,|\hull_\odd^{G-M}(X)|+(\varepsilon-\varepsilon')rk'-2r \\
&\geq|\mathcal{C}_\odd^{G-M}(X)|+\varepsilon'\,|\hull_\odd^{G-M}(X)|
\end{aligned}
\]
which verifies the strengthened Tutte condition in this case.

If $r=1$ and still $|\hull_\odd^G(Y)|\geq k'$, which is the case when $|\hull_\odd^{G-M}(X)|\geq k'$, then the above computation carries through to show $|X|\geq|\mathcal{C}_\odd^{G-M}(X)|+\varepsilon'\,|\hull_\odd^{G-M}(X)|$. Otherwise we must have $\hull_\odd^G(Y)\subseteq B_{k'}(x_1)$. But since $m(x_1)$ was chosen such that $e_{x_1}$ belongs to a matching $N$ covering the ball $B_{k'}(x_1)$, we must have $|X|\geq|\mathcal{C}_\odd^{G-M}(X)|$ since each odd component is matched to an element of $X$ in $N$.

Finally, we explain why $G-M$ has no odd components. Suppose $C\subseteq V(G-M)$ is such a component. As above, list $\{x_1,\ldots,x_r\}$ the set of vertices $x\in A$ such that either $x$ or $m(x)$ is adjacent to $C$. Let $$Y=\{x_1,m(x_1),\ldots,x_r,m(x_r)\}.$$ We have $\mathcal{C}_\odd^G(Y)\ni C$ and $\hull_\odd^G(Y)$ is connected. Then we can rule out the three cases $r=0$, $r\geq2$, and $r=1$ exactly as above working with this $Y\subseteq V(G)$.

\end{proof}

\begin{proof} [Proof of Theorem~\ref{BMTutte}]
Let $G$ be a locally finite Borel graph and $\varepsilon>0$ be such that for every finite $X\subseteq V(G)$ we have $$|X|\geq|\mathcal{C}_\odd(X)|+\varepsilon\,|\hull_\odd(X)|.$$ Let $\varepsilon_0=\varepsilon$ and $k_0=1$. For $n\geq1$ recursively pick $0<\varepsilon_n<\varepsilon_{n-1}$ and $k_n\geq\max\{k_{n-1},\frac{2}{\varepsilon_{n-1}-\varepsilon_n}\}$.

Let $f: \Nat \to \Nat$ be $f(n)=2k_n$. By Lemma~\ref{sparselemma} there exist Borel sets $A_n\subseteq V(G)$ for $n\in\Nat$ such that $\bigcup_nA_n$ is a Borel comeager invariant set and $d_G(x,y)>f(n)$ for distinct vertices $x,y\in A_n$. Let $G_0=G$ and recursively we construct $M_n$ a (partial) matching in $G_n$ covering $A_n\cap V(G_n)$ and let $G_{n+1}=G_n-M_n$. We prove by induction that $G_n$ satisfies $\Tutte_{\varepsilon_n,k_n}$ and that these $M_n$ can be taken to be Borel.

Because $\varepsilon_0=\varepsilon$ and $k_0=1$, we have that $G_0=G$ satisfies $\Tutte_{\varepsilon_0,k_0}$ by assumption. To apply Lemma~\ref{indstep}, for each $x\in A_n\cap V(G_n)$ we let $m(x)$ be the least neighbor of $x$ that can be matched (via $e_x$) to $x$ by a matching covering the ball $B_{k_{n+1}}(x)$ in a fixed Borel linear ordering of $V(G)$. It follows that $M_n=\{e_x:x\in A_n\}$ is Borel. And by Lemma~\ref{indstep}, $G_{n+1}$ satisfies $\Tutte_{\varepsilon_{n+1},k_{n+1}}$.

Thus $M=\bigcup M_n$ is a Borel perfect matching on a comeager invariant subset of $V(G)$.
\end{proof}

Next, we show that non-amenable, vertex transitive graphs satisfy the condition in Theorem~\ref{BMTutte}.

\begin{lemma} \label{NonamenableExpansion}
Let $G$ be an infinite, connected, locally finite, non-amenable, vertex-transitive graph. Then there exists $\ve > 0$ such that for all finite $X \subseteq V(G)$,
\[ |X| \geq |\C_\fin(X)| + \ve \,|\hull_\fin(X)|.\]
In particular, there exists $\ve > 0$ such that for all finite $X \subseteq V(G)$,
\[|X| \geq |\C_\odd(X)| + \ve \,|\hull_\odd(X)|.\]
\end{lemma}
\begin{proof}
Fix a finite set $X \subseteq V(G)$. By Lemma 2.3 of \cite{CL}, the assumption that $G$ is a (connected, infinite) $d$-regular, vertex-transitive graph implies that each finite set of vertices has at least $d$ many edges in its boundary. And so 
\[ \left| E\left(X, \bigcup \C_\fin(X) \right)\right|=\sum_{F\in\mathcal{C}_\text{fin}(X)}\big|E(X, F)\big|\geq d\,|\mathcal{C}_\text{fin}(X)|.\]
Also by the expansion property of non-amenable graphs
\[\left|E\left(X, V(G) \setminus \hull_\fin(X)\right)\right|\geq\delta\,|\hull_\fin(X)|,\]
where $\delta$ is the expansion constant of the graph. Therefore 
\[d\,|X|\geq \left| E\left(X, \bigcup \C_\fin(X) \right) \right| + \Big|E\Big(X, V(G) \setminus \hull_\fin(X)\Big)\Big| \geq d\,|\mathcal{C}_\text{fin}(X)|+\delta\,|\hull_\fin(X)|.\] And so \[|X|\geq |\mathcal{C}_\text{fin}(X)|+\varepsilon\,|\hull_\fin(X)|,\]
where $\varepsilon=\frac{\delta}{d}$.
\end{proof}

As discussed earlier, combining Theorem~\ref{BMTutte} and Lemma~\ref{NonamenableExpansion} immediately yields Theorem~\ref{main}. And Corollary~\ref{maincor} follows because the Schreier graphs of free (Borel) actions of finitely generated non-amenable groups are non-amenable as (Borel) graphs.

\section{Computable Matchings}

Another setting for the study of the definability of matchings in locally finite graphs is computable combinatorics. Qian and Weilacher established connections between the computable and Baire measurable settings for a general class of local labeling problems, including matching.

\begin{definition}
A highly computable graph is a locally finite graph $G$ such that
\begin{itemize}
    \item the vertex set $V(G)\subseteq\Nat$ is computable,
    \item the edge set $E(G)\subseteq\Nat\times\Nat$ is computable, and
    \item the degree function $d:V(G)\to\Nat$ is computable.
\end{itemize}
\end{definition}

We will deal only with infinite graphs $G$, for which we may always assume $V(G)=\Nat$ by a re-coding argument. The condition that the degree function is (finite and) computable is an important restriction that distinguishes this definition from other notions of computable graphs. It is the reason for the word ``highly'' in the term, which we may omit in settings where only graphs of constant finite degree are considered. This condition is equivalent to the condition that the function that maps a vertex to its set of neighbors is computable.

The proof of Theorem~\ref{maincomp} is also via a recursive construction relying on Lemma~\ref{indstep} and will use the same graph-theoretic terms defined in the previous section. First, we prove a computable variant of Tutte's theorem which gives a sufficient condition for a highly computable graph to admit a computable perfect matching.

\begin{theorem}\label{comptutte}
Let $G$ be a highly computable graph on $V(G)=\Nat$ and suppose there exists $\varepsilon>0$ such that for every finite $X\subseteq V(G)$, we have $$|X|\geq|\mathcal{C}_\odd(X)|+\varepsilon|\hull_\odd(X)|.$$ Then $G$ admits a computable perfect matching.
\end{theorem}

\begin{proof}
As before, let $\varepsilon_0=\varepsilon$, which we may assume to be a rational number, and $k_0=1$. Fix computable sequences for $n\geq1$ satisfying $0<\varepsilon_n<\varepsilon_{n-1}$ and $k_n\geq\max\{k_{n-1},\frac{2}{\varepsilon_{n-1}-\varepsilon_n}\}$. Let $f:\Nat\to\Nat$ be $f(n)=2k_n$.

Now we recursively define finite partial matchings $M_n$ covering $\{n\}\cap V(G_n)$ and let $G_{n+1}=G_n-M_n$. We prove by induction that $G_n$ satisfies $\Tutte_{\varepsilon_n,k_n}$ and that these $M_n$ can be computed uniformly in $n$.

Because $\varepsilon_0=\varepsilon$ and $k_0=1$, we have that $G_0=G$ satisfies $\Tutte_{\varepsilon_0,k_0}$ by assumption. If $n\not\in V(G_n)$, then let $M_n=\emptyset$. Otherwise, apply Lemma~\ref{indstep}, where our sparse set is the singleton $\{n\}$. Let $m(n)$ be the least neighbor of $n$ that can be matched (via $e_n$) to $n$ in a matching covering the ball $B_{k_{n+1}}(n)$. It follows that $M_n=\{e_n\}$ is computable uniformly in $n$. And by Lemma~\ref{indstep}, $G_{n+1}$ satisfies $\Tutte_{\varepsilon_{n+1},k_{n+1}}$.

Thus $M=\bigcup M_n$ is a computable perfect matching.
\end{proof}

Theorem~\ref{maincomp} follows because, by Lemma~\ref{NonamenableExpansion}, non-amenable vertex-transitive graphs satisfy the condition in Theorem~\ref{comptutte}. And Corollary~\ref{compcor} follows because the Schreier graphs of free (computable) actions of finitely generated non-amenable groups are non-amenable as (highly computable) graphs.

Similar to what was remarked for Corollary~\ref{maincor}, in the bipartite case Corollary~\ref{compcor} also follows from a previously known variant of Hall's theorem. This effective variant of Hall's theorem is due to Kierstead \cite{effhall}.

\section{Comparison with Bipartite Matching Results}

It is worth comparing Theorem~\ref{BMTutte} with the corresponding result in \cite{MUMatchings}.

\begin{theorem}[Marks-Unger]\label{BMHall}
Let $G$ be a locally finite bipartite Borel graph with bipartition $V(G) = P_0 \sqcup P_1$ (the sets $P_0$ and $P_1$ need not be Borel). Suppose there exists $\ve > 0$ such that whenever $F$ is a finite set contained in either $P_0$ or $P_1$, we have
\[ |N(F)| \geq (1+\ve) |F|.\]
Then $G$ admits a Borel perfect matching on a Borel comeager invariant set.
\end{theorem}

Theorem~\ref{BMTutte} does not immediately imply Theorem~\ref{BMHall}. However, there is still a connection that can be drawn between the two results. Marks and Unger remark that their method also proves the following one-sided matching theorem.

\begin{theorem}[Marks-Unger]\label{BMHall1}
Let $G$ be a locally finite Borel graph on a Polish space $V(G)$ with Borel bipartition $V(G)=P_0\sqcup P_1$, and suppose there exists $\varepsilon>0$ such that for every finite set $X\subseteq P_0$ we have $$|N(X)|\geq(1+\varepsilon)|X|.$$ Then $G$ admits a matching covering $P_0$.
\end{theorem}

Note that a Cantor–Schröder–Bernstein argument shows that Theorem~\ref{BMHall1} implies Theorem~\ref{BMHall} in the special case where the graph has a Borel bipartition. This is the case needed to construct Baire measurable paradoxical decompositions.

We remark that the same method used to prove Theorem~\ref{BMTutte} can be used to prove the following generalization. We omit the proof because it is nearly identical to the proof of Theorem~\ref{BMTutte} and we present no new applications.

\begin{theorem} \label{BMTutte1}
Let $G$ be a locally finite Borel graph on a Polish space $V(G)$, let $P\subseteq V(G)$ be Borel, and suppose there exists $\ve > 0$ such that for every finite set $X \subseteq V(G)$, we have
\[ |X| \geq |\C_\odd^P(X)| + \ve |\hull_\odd^P(X)|.\]
Then $G$ admits a Baire measurable matching that covers $P$.
\end{theorem}

The superscript $P$ indicates that we only consider odd components that are subsets of $P$. If Theorem~\ref{BMTutte} is an analog of Tutte's theorem, then Theorem~\ref{BMTutte1} is an analog of Lovász's theorem \cite{LL} that a set of vertices $P$ in a locally finite graph can be covered by a matching if and only if the removal of finitely many vertices $X$ cannot create more odd components that are subsets of $P$ than there are vertices in $X$. Now Theorem~\ref{BMTutte1} directly implies Theorem~\ref{BMHall1}.

\begin{proof} [Proof of Theorem~\ref{BMHall1} from Theorem~\ref{BMTutte1}]

Let $G$ be a locally finite Borel graph on a Polish space $V(G)$ with Borel bipartition $V(G)=P_0\sqcup P_1$. Let $\varepsilon$ be as in the hypothesis of Theorem~\ref{BMHall1}.

Let $\varepsilon'>0$ be small enough that $$\frac{1+\varepsilon'}{1-\varepsilon'}\leq1+\varepsilon.$$ We show that $G$ satisfies the hypothesis of Theorem~\ref{BMTutte1} for $P=P_0$ and for this value $\varepsilon'$. Let $X\subseteq V(G)$. For the purpose of verifying the inequality, we may assume that $X\subseteq P_1$ since removing elements of $P_0$ only decreases the left side and increases the right side. Because $G$ is bipartite, each odd component of the complement of $X$ contained in $P_0$ is a singleton. Let $Y$ be the union of these singletons. Then $Y\subseteq P_0$, with $|\C_\odd^P(X)|=|Y|$ and $|\hull_\odd^P(X)|=|X|+|Y|$. We therefore have, by assumption, $|X|=|N(Y)|\geq(1+\varepsilon)|Y|$. Putting this all together we see that
\begin{align*}
|\C_\odd^P(X)| + \ve' \,|\hull_\odd^P(X)|&=|Y|+\varepsilon'(|X|+|Y|)\\
&=(1+\varepsilon')\,|Y|+\varepsilon'\,|X|\\
&\leq(1+\varepsilon')\,\frac{|X|}{1+\varepsilon}+\varepsilon'\,|X| \\
&\leq(1-\varepsilon')\,|X|+\varepsilon'\,|X|\\
&=|X|
\end{align*}
which verifies the hypothesis of Theorem~\ref{BMTutte1}. And so $G$ has a Baire measurable matching covering $P_0$. This completes the argument.

\end{proof}

Therefore, in some sense, the Baire measurable expansive Tutte results do generalize the Baire measurable expansive Hall results.

\section{Balanced orientations}

The proof of Theorem~\ref{orientation} is quite straightforward and is an adaptation of the ideas in Section 5 of \cite{BHTOrientations}. Given any graph $G$ with only even degrees, we define an auxiliary bipartite graph $G^*$ such that perfect matchings of $G^*$ induce balanced orientations of $G$. The following definition is taken essentially verbatim from \cite{BHTOrientations} and is repeated here for the convenience of the reader.
\begin{definition}
Let $G$ be a graph with only even degrees. The graph $G^*$ has a vertex for every edge $e$ of $G$ and $\deg(v)/2$ vertices for every vertex $v$ of $G$, i.e.
\[ V(G^*) = \big\{x_e: e \in E(G)\big\} \cup \Big\{v_i: v \in V(G), i \in \big[\deg(v)/2\big]\Big\}.\]
Then every vertex corresponding to a former edge is joined to all copies of its former endpoints:
\[ E(G^*) = \Big\{x_{uv} v_i: u v \in E(G), i \in \big[\deg(v)/2\big]\Big\}.\]
\end{definition}
Observe that any perfect matching of $G^*$ induces a balanced orientation of $G$ by orienting an edge $e \in E(G)$ toward its endpoint $v$ if and only if $x_e$ and $v_i$ are matched in $G^*$ for some $i \in [\deg(v)/2]$. In the case when $G$ is a Borel graph, it is straightforward to put an appropriate Polish topology on $V(G^*)$ such that Baire measurable perfect matchings of $G^*$ yield Baire measurable balanced orientations of $G$. In order to show that the Borel bipartite graph $G^*$ has a Baire measurable perfect matching, we will apply Theorem~\ref{BMHall} of \cite{MUMatchings} (which is an analogue of our Theorem~\ref{main} in the bipartite setting).

\begin{proof}[Proof of Theorem~\ref{orientation}]
The proof is an adaptation of the proof of Lemma 25 from \cite{BHTOrientations}. Write $\pi: V(G^*) \to V(G) \cup E(G)$ for the projection function. Let $\delta > 0$ be the expansion constant for the non-amenable graph $G$, and let $d$ be a bound on the degrees. Suppose that $F \subseteq V(G^*)$ is a finite set of vertex-type vertices. Then
\[ \begin{aligned}
|N_{G^*}(F)| &= \frac{1}{2} \sum_{u \in \pi(F)} \deg(u) + \frac{1}{2} \,|E(\pi(F), V(G) \setminus \pi(F))| \\
&\geq \sum_{u \in \pi(F)} \frac{\deg(u)}{2} + \frac{\delta}{2}\,|\pi(F)| \\
&\geq |F| + \frac{\delta}{d}\,|F|.
\end{aligned}\]
Now suppose that $F \subseteq V(G^*)$ is a finite set of edge-type vertices (that is, $F \subseteq E(G)$), and let $S$ denote the set of vertices $u \in V(G)$ which are incident to some edge $e \in F$. Then
\[ \begin{aligned}
|N_{G^*}(F) &= \sum_{u \in S} |\pi^{-1}(u)| \\
&= \sum_{u \in S} \frac{\deg(u)}{2} \\
&= |E(S, S)| + \frac{1}{2} \,|E(S, V(G) \setminus S)| \\
&\geq |F| + \frac{\delta}{2} \,|S| \\
&\geq |F| + \frac{\delta}{2d} \,|F|.
\end{aligned}\]
If we choose $\ve > 0$ such that $\ve < \frac{\delta}{2d}$, then the hypotheses for Theorem~\ref{BMHall} hold. So $G^*$ has a Baire measurable perfect matching, and this implies that $G$ has a Baire measurable balanced orientation.
\end{proof}

\end{document}